\documentclass[preprint,12pt]{elsarticle} 
\usepackage[english]{babel}
 
\usepackage{enumerate}            
\usepackage{amsmath}  
\usepackage{amssymb}    
\usepackage{graphicx}       
\usepackage{epsfig}            

\usepackage{caption}
\usepackage{float}
\usepackage{tabulary} 

\usepackage{enumerate}     
\usepackage{amsmath}

\journal{journal}  
\begin{document}

\makeatletter
\newcommand{\norm}[1]{\left\lVert#1\right\rVert}
\makeatother
  

\begin{abstract}
The conceptually new approach based on the logarithmic norm to design of robust adaptive state-feedback controller for linear time-varying (LTV) systems under system's modeling uncertainty and nonlinear external disturbance is proposed. This controller, consisting of two independent parts - adaptive and robust ones - globally asymptotically stabilizes every LTV system regardless how large the disturbance is. 
\end{abstract} 

\title{Robust adaptive state-feedback control of linear time-varying systems under both potentially unbounded system's modeling uncertainty and external disturbance} 
\author{R.~Vrabel}
\ead{robert.vrabel@stuba.sk}
\address{Slovak University of Technology in Bratislava, Institute of Applied Informatics, Automation and Mechatronics,  Bottova 25,  917 01 Trnava,   Slovakia}

\newtheorem{thm}{Theorem}
\newtheorem{lem}[thm]{Lemma}
\newtheorem{defi}[thm]{Definition}
\newtheorem{cor}[thm]{Corollary}
\newdefinition{rmk}{Remark}
\newdefinition{ex}{Example}
\newproof{pf}{Proof}
\newproof{pot1}{Proof of Theorem \ref{thm1}}
\newproof{pot2}{Proof of Theorem \ref{thm2}}

\pagestyle{headings}

\maketitle

\section[Introduction]{Introduction and preliminary results} 

The unknown and/or unmeasurable variations of the process parameters and also the external disturbances acting upon the controlled variables can degrade the performances of any control system. Adaptive control strategies cover a set of techniques providing a systematic approach for automatic adjustment of controllers in real time, in order to achieve or to maintain a desired level of control system performance when the parameters of the plant dynamic model are initially uncertain and/or change in time.  Robust control differs from adaptive control, while adaptive control is concerned with control law changing itself, robust control guarantees that if the  (potentially unknown) disturbances are within given constraints, bounds or eventually asymptotics as $t\to\infty,$ the control law need not be changed. Trajectory tracking is one of the important, if not the most important, motion control problems, which not only requires a designed controller \cite{Worthmann} but also has to robustly stabilize the nonlinear system against the system's modeling uncertainties and external disturbances \cite{Ce}. The stable tracking of the desired state trajectory $x_d(t)$ can be transformed into analysis of tracking error dynamics $e(t)=x(t)-x_d(t)$ and its stabilizability at the equilibrium $e=0$ \cite{Vrabel}, \cite{Xin}, \cite{Zhou}. Thus, without loss of generality, as the desired value we will consider here the origin $x=0$ of $\mathbb{R}^n.$

A renewed interest on robust control appeared in the late 1970s and early 1980s \cite{Safonov} in the connection with the problem of plant uncertainty \cite{Dorato} and soon developed a various methods for dealing with bounded uncertainties \cite{Ioannou}. Among them, for example, the sliding mode control \cite{Feng}, \cite{Golestani}, \cite{Levant}, \cite{Utkin}, \cite{Xin}, \cite{Zhang},  the  backstepping method  \cite{Kanellakopoulos}, \cite{Wu}, finite-time control technique \cite{Bhat}, \cite{Ce}, \cite{Golestani} and the fuzzy control \cite{Mao}, \cite{Zhou}.
However, none of these techniques are well suited for the systems with {\bf unbounded uncertainties} - system's modeling uncertainty and external disturbance. Therefore, it is not currently possible to find scientific sources to compare the results achieved by different methods and approaches.

The feedback control is basically used in the conventional control systems to reject the effect of all these mentioned disturbances upon the controlled variables and to bring the plants back to their desired values, $x=0.$ Combining and integrating above-mentioned notions the robust adaptive state-feedback control is designed in Corollary~\ref{corollary1}, which is 
\begin{itemize}
\item[i)] adaptive with regard to the system's modeling uncertainty $\Delta(t)$ and independent of it if (possibly unbounded) uncertainty satisfies the constraint given in Theorem~\ref{main_theorem_disturbances}-Assumption~A1, 
\item[ii)] robust with regard to external disturbance $\omega(x,t)$ satisfying Assumption~A3 of Theorem~\ref{main_theorem_disturbances}, and
\item[iii)] these two components are independent of each other.
\end{itemize}
Specifically, we will consider the state-feedback control system
\begin{equation} \label{basic_model}
\dot x=\big[A(t)+\Delta(t)\big]x+Bu+\omega(x,t),\  u=K(t)x,\ t\geq t_0,
\end{equation}
where $A,$ $\Delta,$ and $K:$ $[t_0,\infty)\to\mathbb{R}^{n\times n}$ are (piecewise) continuous matrix functions, representing known part of system's dynamics, system's modeling uncertainty, and state-feedback gain matrix, respectively, $B$ is a constant matrix ($n\times n$) and the (potentially unknown) disturbance $\omega:$ $\mathbb{R}^n\times[t_0,\infty)\to\mathbb{R}^n$ is (piecewise) continuous for all $x\in\mathbb{R}^n$ and all $t\geq t_0.$ We do not assume here that $\omega(0,t)\equiv0,$ that is, $x=0$ may not be the equilibrium point for dynamical model (\ref{basic_model}). This is reason, why the approach based on the Lyapunov function is not suitable, and it is useful to find other ways to find out when the system (\ref{basic_model}) solutions converge to zero as $t\to\infty$. 

The first main goal of the present paper is to establish (in Theorem~\ref{main1}), in terms of the logarithmic norm, the sufficient conditions to be the system with system's modeling uncertainty but without external disturbance,
\begin{equation} \label{basic_model_without}
\dot x=\big[A(t)+\Delta(t)\big]x+Bu,\ u=K(t)x,\ t\geq t_0,
\end{equation}
(uniformly, asymptotically, uniformly asymptotically) stable in the sense of classical definitions \cite[p.~149]{Khalil}. These various types of stability can be expressed in the terms of a fundamental matrix \cite[p.~54]{Coppel2}.
\begin{thm}\label{thm:stabilityLTV}
Let $\Phi(t)$ be a fundamental matrix  for $\dot x={\cal A}(t)x,$ $t\geq t_0.$ Then the system $\dot x={\cal A}(t)x$ is 
\begin{itemize} 
\item[(S)] stable if and only if there exists a positive constant $M$ such that
\[
\norm{\Phi(t)}\leq M \ \mathrm{for\ all}\ t\geq t_0,
\]
\item[(US)] uniformly stable if and only if there exists a positive constant $M$ such that
\begin{equation*}
\norm{\Phi(t)\Phi^{-1}(\tau)}\leq M\ \mathrm{for}\ t_0\leq\tau\leq t<\infty,
\end{equation*}
\item[(AS)] asymptotically stable if and only if
\begin{equation*}
\norm{\Phi(t)}\to 0\ \mathrm{as}\ t\to\infty,
\end{equation*}
\item[(UAS)]  uniformly asymptotically stable ($\Leftrightarrow$ uniformly exponentially stable) if and only if there exist positive constants
$M$, $\alpha$ such that
\begin{equation*}
\norm{\Phi(t)\Phi^{-1}(\tau)}\leq Me^{-\alpha(t-\tau)}\ \mathrm{for}\ t_0\leq\tau\leq t<\infty.
\end{equation*}
\end{itemize}
\end{thm}
The another goal of this paper is at the same time to determine the sufficient conditions ensuring the convergence of all solutions of the system with system's modeling uncertainty and external disturbance, (\ref{basic_model}), to $0$ as $t\to\infty,$ in Theorem~\ref{main_theorem_disturbances} and Corollary~\ref{corollary1}.

We will derive our results for unspecified vector norm on $\mathbb{R}^n$, $\norm{\cdot}.$ For the matrices as an operator norm the induced norm will be used, $\norm{{\cal A}}=\max\limits_{\norm{x}=1}\norm{{\cal A}x}$.  We use for both vector norm and matrix operator norm the same notation but it will always be clear from the context that norm is being used. In particular cases we will consider the three most common vector norm: 
\begin{equation}\label{norms}
\norm{x}_1=\sum\limits_{i=1}^{n}|x_i|, \qquad \norm{x}_2=\left(\sum\limits_{i=1}^{n}x_i^2\right)^{1/2}, \qquad \norm{x}_{\infty}=\max_{1 \leq i \leq n}{|x_i|}.
\end{equation}

By $\mu[{\cal A}(t)],$ $t\geq t_0,$ we denote the logarithmic norm of an $n\times n$ matrix ${\cal A}(t)$ defined as
\[
\mu[{\cal A}(t)]\triangleq\lim\limits_{h\to 0^+}\frac{\norm{I_n+h{\cal A}(t)}-1}{h},
\]
where $I_n$ is the identity on $\mathbb{R}^n.$ The logarithmic norm is not a norm in usual sense. While the matrix norm $\norm{{\cal A}}$ is always positive if ${\cal A}\neq0,$  the logarithmic norm $\mu[{\cal A}]$ may also take negative values, e.~g. for the Euclidean vector norm $\norm{\cdot}_2$ and when ${\cal A}$ is negative definite because $\frac12({\cal A}+{\cal A}^T)$ is also negative definite, 
\cite[Corollary~14.2.7.]{Harville} and Lemma~\ref{lognorm:properties}. Therefore, the logarithmic norm does not satisfy the axioms of a norm. 

The fundamental advantage of approach based on the use of logarithmic norm is the fact that to estimate the norm of transition matrix $\Phi(t)\Phi^{-1}(\tau)$ for LTV system $\dot x={\cal A}(t)x$  we do not need to know the fundamental matrix explicitly. Moreover, because $\mu[{\cal A}(t)]\leq\norm{{\cal A}(t)}$ and $\mu[{\cal A}(t)]$ can take also negative values, we can obtain the stronger results as those achieved in the now classical results (and their numerous variations) regarding persistence of the stability properties of perturbed system, $\dot x=\big[{\cal A}(t)+{\cal B}(t)\big]x,$ where it is assumed that the perturbing term satisfies
\begin{itemize}
\item[] $\int\limits_{t_0}^\infty\norm{{\cal B}(s)}ds<\infty$ \cite[p.~65]{Coppel2}, \cite[p.~133]{Rugh}, or
\item[] $\norm{{\cal B}(t)}\to 0$ as $t\to\infty$  \cite[p.~354]{Khalil}, \cite[p.~70]{Coppel2}.
\end{itemize}
Note, that none of these conditions are required in this paper.

\centerline{}
 
In Table~\ref{table:norms} and Lemma~\ref{lognorm_properties} we summarize properties of the logarithmic norm useful for the stability analysis of dynamical systems.

\centerline{}

\begin{minipage}{\textwidth}
\begin{lem}[\cite{Afanasiev}, \cite{Desoer_Vidyasagar}] \label{lognorm:properties} 
For the norms (\ref{norms}) we have:
\begin{table}[H]
\caption{Logarithmic norms for the vector norms $\norm{\cdot}_1,$ $\norm{\cdot}_2$ and $\norm{\cdot}_\infty$}
\centering 
\scalebox{0.9}{
\begin{tabulary}{\linewidth}{CC}

\hline\hline                        
Norm of the vector ($\norm{x}_i$) & Logarithmic norm ($\mu_i[{\cal A}]$)  \\  
\hline               
$\norm{x}_1=\sum\limits_{i=1}^{n}|x_i|$  & $\mu_1[{\cal A}]=\max\limits_{1 \leq j \leq n}\left(a_{jj}+\sum\limits_{i\neq j}|a_{ij}|\right)$   \\
 &  \\
$\norm{x}_2=\left(\sum\limits_{i=1}^{n}x_i^2\right)^{1/2}$ & $\mu_2[{\cal A}]=\frac12\lambda_{\max}\left({{\cal A}+{\cal A}^T}\right)$   \\
 &  \\
$\norm{x}_{\infty}=\max\limits_{1 \leq i \leq n}{|x_i|}$  & $\mu_{\infty}[{\cal A}]=\max\limits_{1 \leq i \leq n}\left(a_{ii}+\sum\limits_{j\neq i}|a_{ij}|\right)$   \\     
\hline
\end{tabulary}
}
\label{table:norms}
\end{table}
\end{lem} 
\end{minipage}

\centerline{}

\centerline{}

In Table~\ref{table:norms} and elsewhere in the paper, the superscript 'T' denotes transposition,  the number $\lambda_{\max}({\cal A}+{\cal A}^T)$ is the maximum eigenvalue of the matrix ${\cal A}+{\cal A}^T.$

\begin{lem}[\cite{Desoer_Vidyasagar, Desoer2, Soderlind1, Soderlind2}]\label{lognorm_properties}
For any given $n\times n$ real matrix ${\cal A}$ and ${\cal B},$
\begin{itemize}
\item[P1)] $-\norm{{\cal A}}\leq-\mu[-{\cal A}]\leq\mu[{\cal A}]\leq\norm{{\cal A}};$
\item[P2)]  $\mu[{\cal A}+{\cal B}]\leq\mu[{\cal A}]+\mu[{\cal B}]$ and $\vert\mu[{\cal A}]-\mu[{\cal B}]\vert\leq\norm{{\cal A}-{\cal B}};$
\item[P3)] \cite{Desoer2}  Let $\Phi(t)$ is a fundamental matrix for $\dot x={\cal A}(t)x.$ Then 
\[
e^{-\int\limits_{\tau}^t \mu[-{\cal A}(s)]ds}\leq\norm{\Phi(t)\Phi^{-1}(\tau)}\leq e^{\int\limits_{\tau}^t \mu[{\cal A}(s)]ds}
\]
for all $t_0\leq\tau\leq t<\infty;$
\item[P4)]\cite[p.~34]{Desoer_Vidyasagar} The solution of $\dot x={\cal A}(t)x$ satisfies for all $t\geq t_0$ the inequalities
\[
\norm{x(t_0)}e^{-\int\limits_{t_0}^t \mu[-{\cal A}(s)]ds}\leq\norm{x(t)}\leq\norm{x(t_0)}e^{\int\limits_{t_0}^t \mu[{\cal A}(s)]ds}.
\]
\end{itemize}
\end{lem}
\begin{rmk}
Property~P3 allows to estimate the norm of the state-transition matrix without knowing the fundamental matrix, purely on the basis of the matrix ${\cal A}(t)$ entries, which can be especially useful if ${\cal A}(t)$ is a non constant matrix. Moreover, because the logarithmic norm $\mu[{\cal A}(t)]$ can attain also negative values, the estimations above provide information about the actual growth or decay rate in the system. Property~P2 (its second part) implies that $\mu[{\cal A}(t)]$ is (piecewise) continuous function on $[t_0,\infty)$ if such is also the matrix function ${\cal A}:\, [t_0,\infty)\to\mathbb{R}^{n\times n}.$ 
\end{rmk} 

The following example shows that the value $\mu[{\cal A}]$ may depend on the used vector norm.
\begin{ex}\label{example_norms} \cite[p.~56]{Afanasiev} If
\begin{itemize}
\item[a)]  
\[
{\cal A}_1=\begin{bmatrix}
-11 & 10 \\
 2 & -3 
\end{bmatrix},
\]
then $\mu_1[{\cal A}_1]= 7,$ $\mu_2[{\cal A}_1]=0.211$ and $\mu_\infty[{\cal A}_1]=-1.$
\item[b)]
\[
{\cal A}_2=\begin{bmatrix}
-11 & 2 \\
 10 & -3 
\end{bmatrix},
\]
then $\mu_1[{\cal A}_2]= -1,$ $\mu_2[{\cal A}_2]=0.211$ and $\mu_\infty[{\cal A}_2]=7.$
\item[c)]  
\[
{\cal A}_3=\begin{bmatrix}
-1 & 3 \\
 -3 & -2 
\end{bmatrix},
\]
then $\mu_1[{\cal A}_3]= 2,$ $\mu_2[{\cal A}_3]=-1$ and $\mu_\infty[{\cal A}_3]=2.$
\end{itemize}
Thus, we can verify whether the LTI system $\dot x={\cal A}_ix, $ $i=1,2,3$ is asymptotically stable or not by means of the vector norm with negative value of $\mu[{\cal A}_i].$ It is worth to note that for any Hurwitz matrix ${\cal A},$ using a vector norm $\norm{x}_H=(x^THx)^{1/2},$ where the symmetric positive definite matrix $H$ satisfies the Lyapunov equation ${\cal A}^TH + H{\cal A} = -2I_n$, the corresponding $\mu_H[{\cal A}]=-1/\lambda_{\max}(H),$  for details see e.~g. \cite{Hu_Hu}  and  \cite{Hu_Liu} (Lemma~2.3). The stability analysis based on the logarithmic norm becomes a topological notion unlike the spectrum of matrices which is topological invariant.
\end{ex} 

\section{Main results}

We have the following result regarding different types of stability.
\begin{thm}\label{main1}
Consider the control system with system's modeling uncertainty
\[
\dot x=\big[A(t)+\Delta(t)\big]x+Bu+\omega(x,t),\  u=K(t)x,\ t\geq t_0. 
\]
Suppose that for some vector norm in $\mathbb{R}^n$
\[
\int\limits_{t_0}^\infty |\mu[\Delta(s)]|ds<\infty.
\]
Then the control system (\ref{basic_model_without}) is
\begin{itemize}
\item stable if $\limsup\limits_{t\to\infty}\int\limits_{t_0}^t\mu[A(s)+BK(s)]ds<\infty;$
\item uniformly stable if $\mu[A(t)+BK(t)]\leq 0$ for $t\geq t_0;$ 
\item (globally) asymptotically stable if $\int\limits_{t_0}^\infty\mu[A(s)+BK(s)]ds=-\infty;$
\item (globally) uniformly asymptotically stable if $\mu[A(t)+BK(t)]\leq -\alpha<0$ for $t\geq t_0;$ 
\item unstable if $\liminf\limits_{t\to\infty}\int\limits_{t_0}^t\mu[-A(s)-\Delta(s)-BK(s)]ds=-\infty.$
\end{itemize}
\end{thm}
\begin{pf}
The solutions of (\ref{basic_model_without}) can be expressed as
\begin{equation*}
x(t)=\Phi(t)\Phi^{-1}(t_0)x(t_0),\ \ t\geq t_0,
\end{equation*}
where $\Phi(t)$ is a fundamental matrix for LTV system (\ref{basic_model_without}), and estimated for all $t\geq t_0$ as
\[
\norm{x(t)}\leq\norm{\Phi(t)\Phi^{-1}(t_0)}\norm{x(t_0)}\leq\norm{x(t_0)}e^{\int\limits_{t_0}^t \mu[A(s)+\Delta(s)+BK(s)]ds}
\]
\[
\leq\norm{x(t_0)}e^{\int\limits_{t_0}^t |\mu[\Delta(s)]|ds}e^{\int\limits_{t_0}^t \mu[A(s)+BK(s)]ds}.
\]
The rest of proof follows now immediately by using Theorem~\ref{thm:stabilityLTV} and Lemma~\ref{lognorm_properties}.
\end{pf}

\begin{thm}\label{main_theorem_disturbances}
Consider the state-feedback control system with the system's modeling uncertainty and external disturbance
\[
\dot x=\big[A(t)+\Delta(t)\big]x+Bu+\omega(x,t),\  u=K(t)x,\ t\geq t_0.
\]
Let for some vector norm in $\mathbb{R}^n$ and state-feedback gain matrix $K(t),$
\begin{itemize}
\item[A1)] $\int\limits_{t_0}^\infty |\mu[\Delta(s)]|ds<\infty;$
\item[A2)]  $\mu[A(t)+BK(t)]<0$ in some left neighborhood of $\infty;$ and
\item[A3)] for all $x\in\mathbb{R}^n$ and all $t\geq t_0$ is $\norm{\omega(x,t)}\leq\norm{\tilde{\omega}(t)},$  
\[
\lim\limits_{t\to\infty}\frac{\norm{\tilde{\omega}(t)}}{\mu[A(t)+BK(t)]}=0,
\]
that is, $\norm{\tilde{\omega}(t)}=o\left(\mu[A(t)+BK(t)]\right)$ as $t\to\infty.$
\end{itemize}
Then for the state-feedback control law $B(t)x$ satisfying 
\begin{itemize}
\item[A4)] $\int\limits_{t_0}^\infty \mu[A(s)+BK(s)]ds=-\infty,$
\end{itemize}
all solution of (\ref{basic_model}) converge to $0$ as $t\to\infty.$
\end{thm}
\begin{pf}
By the variation of constants formula we have for any solution  of (\ref{basic_model}) 
\begin{equation*}
x(t)=\Phi(t)\Phi^{-1}(t_0)x(t_0)+\Phi(t)\int\limits_{t_0}^t \Phi^{-1}(\tau)\omega(x(\tau),\tau)d\tau,
\end{equation*}
that is, according to Lemma~\ref{lognorm_properties} and Assumption~A3,
\[
\norm{x(t)}\leq \norm{x(t_0)}e^{\int\limits_{t_0}^t\mu[\Delta(s)]ds}e^{\int\limits_{t_0}^t\mu[A(s)+BK(s)]ds}
\]
\[
+\int\limits_{t_0}^t e^{\int\limits_{\tau}^t\mu[\Delta(s)]ds}e^{\int\limits_{\tau}^t\mu[A(s)+BK(s)]ds}\norm{\tilde \omega(\tau)}d\tau.
\]
Obviously, by the Assumptions~A1 and~A4, 
\[
\norm{x(t_0)}e^{\int\limits_{t_0}^t\mu[\Delta(s)]ds}e^{\int\limits_{t_0}^t\mu[A(s)+BK(s)]ds}\to 0\ \mathrm{as}\ t\to\infty
\]
for an arbitrary $x(t_0),$ proving the (global) asymptotic stability of the equilibrium $x=0$ of the system (\ref{basic_model_without}). Because of the absolute convergence of $\int\limits_{t_0}^\infty \mu[\Delta(s)]ds$ (Assumption~A1), $e^{\int\limits_{\tau}^t\mu[\Delta(s)]ds}$ is uniformly bounded for all $t\geq\tau\geq t_0$ and so it remains to analyze the second term on the left-hand side of the above inequality. We have  
\[
\int\limits_{t_0}^t e^{\int\limits_{\tau}^t\mu[A(s)+BK(s)]ds}\norm{\tilde \omega(\tau)}d\tau
\]
\[
= e^{\int\limits_{t_0}^t\mu[A(s)+BK(s)]ds}\int\limits_{t_0}^t e^{-\int\limits_{t_0}^{\tau}\mu[A(s)+BK(s)]ds}\norm{\tilde \omega(\tau)}d\tau
\]
\[
=\frac{\int\limits_{t_0}^t e^{-\int\limits_{t_0}^{\tau}\mu[A(s)+BK(s)]ds}\norm{\tilde \omega(\tau)}d\tau}{e^{-\int\limits_{t_0}^t\mu[A(s)+BK(s)]ds}},
\]
and the L'Hospital rule, allowed by Assumptions~A2 and A4, yields
\[
\lim\limits_{t\to\infty}\frac{\frac{d}{dt}\int\limits_{t_0}^t e^{-\int\limits_{t_0}^{\tau}\mu[A(s)+BK(s)]ds}\norm{\tilde \omega(\tau)}d\tau}{\frac{d}{dt}e^{-\int\limits_{t_0}^t\mu[A(s)+BK(s)]ds}}
\]
\[
=\lim\limits_{t\to\infty}\frac{e^{-\int\limits_{t_0}^t\mu[A(s)+BK(s)]ds}\norm{\tilde w(t)}}{e^{-\int\limits_{t_0}^t\mu[A(s)+BK(s)]ds}(-\mu[A(t)+BK(t)])}=-\lim\limits_{t\to\infty}\frac{\norm{\tilde \omega(t)}}{\mu[A(t)+BK(t)]}.
\]
This, together with Assumption~A3, gives the claim of Theorem~\ref{main_theorem_disturbances}.
\end{pf}

In the following section it is shown that the robust adaptive state-feedback controller $B(t)x$ consists of two independent parts, namely, the adaptive part for the system's modeling uncertainty and the robust part for the external disturbance, reflecting the different nature of the adaptive and robust control.

\section[Construction of the controller]{Construction of the robust adaptive state-feedback controller in the Euclidean norm}

In this section we show the use of Theorem~\ref{main_theorem_disturbances} for the most frequently used norm in the state space $\mathbb{R}^n$, namely, the Euclidean norm. Of course, depending on the particular system, it may be more appropriate to use a different vector norm as indicated and justified in Example~\ref{example_norms} to adapt the choice of gain matrix $K(t).$

\begin{cor}\label{corollary1}
For the Euclidean vector norm $\norm{\cdot}_2$ and for every disturbance $\omega(x,t)$ satisfying $\norm{\omega(x,t)}_2\leq\norm{\tilde{\omega}(t)}_2$ for all $(x,t)\in\mathbb{R}^n\times[t_0,\infty),$ for every system's modeling uncertainty $\Delta(t)$ satisfying Assumption~A1 of Theorem~\ref{main_theorem_disturbances} and invertible control matrix $B$ there exists a robust adaptive state-feedback control law $u=K(t)x$ such that all solutions of closed-loop system
\[
\dot x=\big[A(t)+\Delta(t)\big]x+Bu+\omega(x,t),\ t\geq t_0,
\]
converge to $0$ as $t\to\infty.$
\end{cor}
\begin{pf}
At first, let us decompose the system matrix $A(t)$ into its symmetric and skew-symmetric part, 
\[
A(t)=\frac12(A(t)+A^T(t))+\frac12(A(t)-A^T(t))\triangleq A_{\mathrm{sym}}(t)+A_{\mathrm{skew}}(t),
\]
and define the gain matrix $K(t)$ as 
\[
K(t)=B^{-1}\left(\underbrace{-A_{\mathrm{sym}}(t)+\mathrm{diag}\{\lambda_1,\dots,\lambda_n\}}_{\mathrm{adaptive\ part\ of\ control}}+\underbrace{\mathrm{diag}\{\gamma_1(t),\dots,\gamma_n(t)\}}_{\mathrm{robust\ part\ of\ control}}\right),
\]
where 
\begin{itemize}
\item[c1)] $\mathrm{diag}\{\lambda_1,\dots,\lambda_n\}$ and $\mathrm{diag}\{\gamma_1(t),\dots,\gamma_n(t)\}$ denote a diagonal square matrix with $\lambda_i$'s and $\gamma_i$'s on the main diagonal and the entries outside the main diagonal are all $0,$ respectively;
\item[c2)] all $\lambda_i$ are negative real numbers and the (piecewise) continuous functions $\gamma_i(t)$ defined on $[t_0,\infty)$ are chosen such that $\norm{\tilde{\omega}(t)}_2=o(\gamma_i(t))$ as $t\to\infty$ for all $i=1,\dots,n$ and 
\item[c3)] $\int\limits_{t_0}^\infty\Gamma(s)ds=-\infty$ for $\Gamma(t)=\max\{ \lambda_i+\gamma_i(t);\, i=1,\dots,n\}.$ 
\end{itemize}
Now, taking into account that $\mu_2[A(t)+BK(t)]=\Gamma(t),$  the claim of Corollary~\ref{corollary1} follows by Theorem~\ref{main_theorem_disturbances}.
\end{pf}

\subsection{Simulation experiment}
\begin{figure}[ht] 
\captionsetup{singlelinecheck=off}
   \centerline{
    \hbox{
     \psfig{file=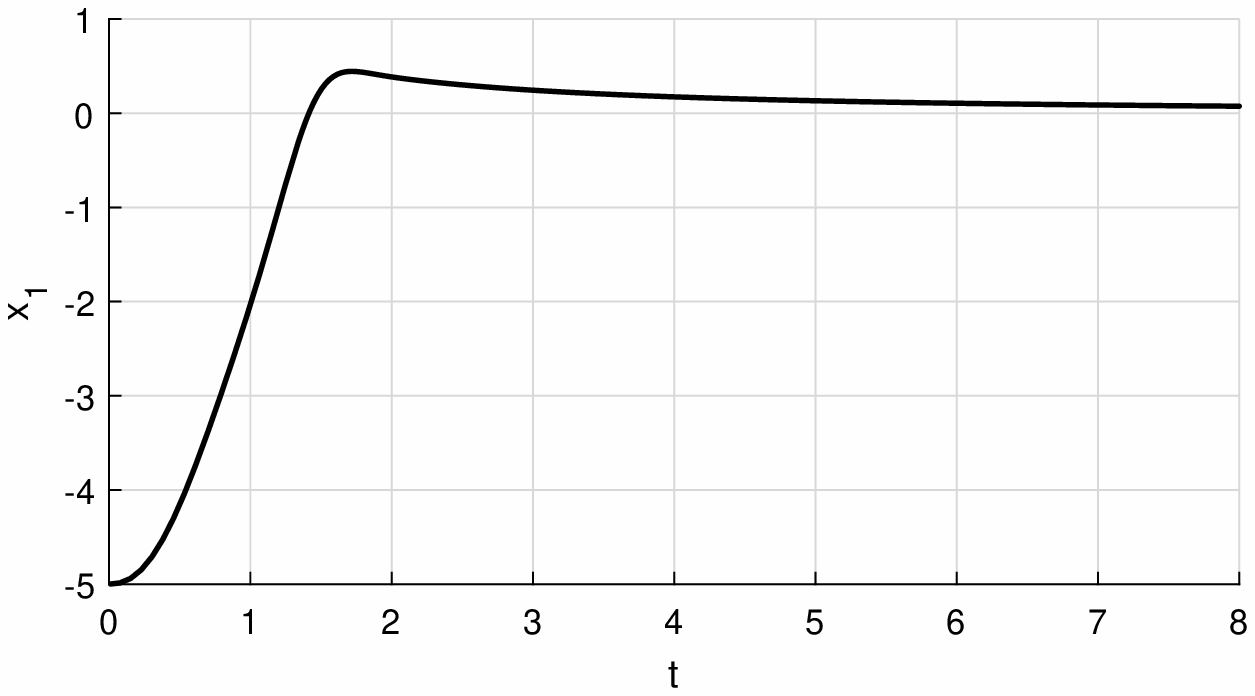,width=5.0cm, clip=}
     \hspace{1.cm}
     \psfig{file=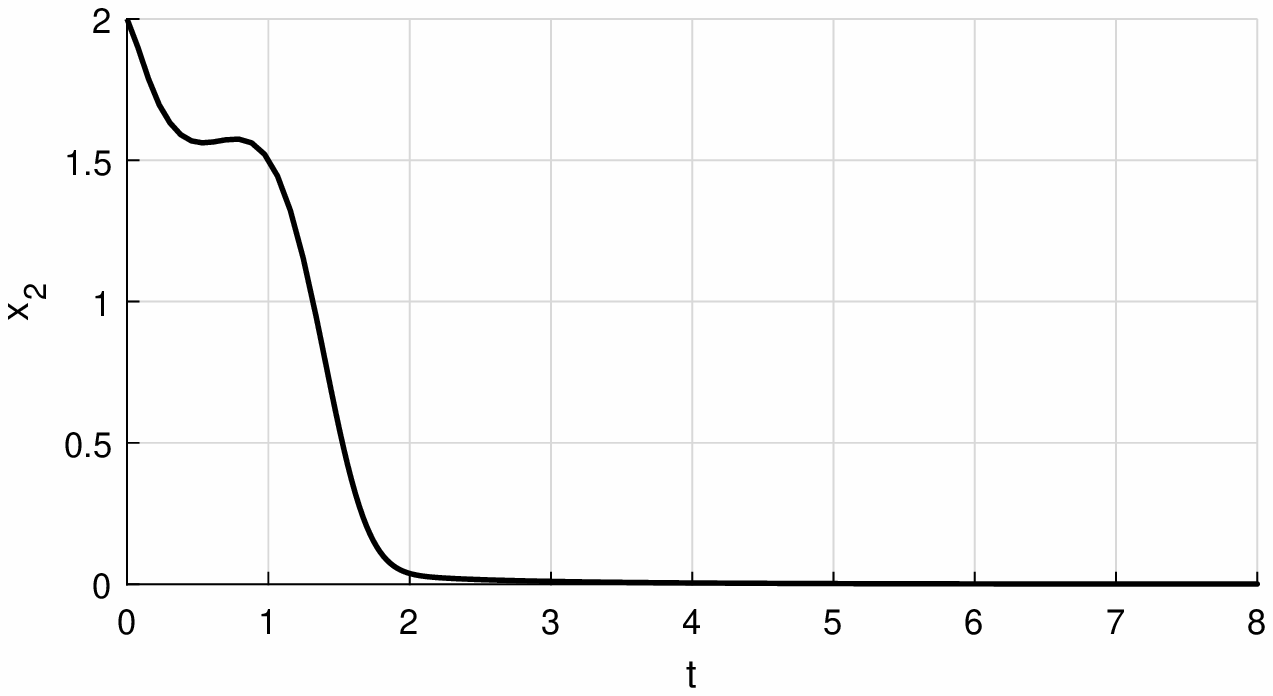,width=5.0cm,clip=}
    }
   }
   
\caption{The simulation result for the system (\ref{example1}) with $\omega(x,t)=(t^{11/4}\cos x_1,\, 1)^T $ and the initial state $x(0)=(-5,\, 2)^T$ demonstrating the effectiveness of proposed controller.}
\label{simulation_x_1_2}
\end{figure}

\begin{ex}\label{example_all_to_zero}
To illustrate the theory developed above let us consider the control system
\begin{equation}\label{example1}
\dot x=\big[A(t)+\Delta(t)\big]x+Bu+\omega(x,t),\  u=K(t)x,\ t\geq 0,
\end{equation}
where
\[
A(t)=\begin{bmatrix}
t & \sin t \\
 t^{1/2} & 1 
\end{bmatrix},\
\Delta(t)=
\begin{bmatrix}
1/(1+t^2) & t \\
 -t & 0 
\end{bmatrix},\ B=I_2\ \mathrm{and}\ \omega(x,t)=
\begin{bmatrix}
o(t^3)  \\
 1 
\end{bmatrix},
\]
that is, the system with unbounded on the time interval  $[0,\infty)$ system's modeling uncertainty $\Delta$ and unknown unbounded external disturbance $\omega.$  

Obviously, 
\[
\mu_2[\Delta(t)]=\frac12\lambda_{\max}\left({\Delta(t)+{\Delta(t)}^T}\right)=1/(1+t^2)
\]
and so
\[
 \int\limits_0^\infty |\mu_2[\Delta(s)]|ds=\pi/2<\infty
\]
(Assumption~A1 of Theorem~\ref{main_theorem_disturbances}). The only information we have about the external disturbance is that $\norm{\omega(x,t)}_2\leq (o^2(t^3)+1)^{1/2}\ (=\norm{\tilde \omega(t)}_2).$ 

By Corollary~\ref{corollary1}, we can define the control matrix as follows:
\[
K(t)=-A_{\mathrm{sym}}(t)+\mathrm{diag}\{\lambda_1,\, \lambda_2\}+\mathrm{diag}\{\gamma_1(t),\, \gamma_2(t)\}
\]
\[
=
-\begin{bmatrix}
t & (t^{1/2}+\sin t)/2 \\
 (t^{1/2}+\sin t)/2 & 1 
\end{bmatrix}+
\begin{bmatrix}
-1& 0 \\
 0 & -1 
\end{bmatrix}
\]
\[
+
\begin{bmatrix}
-t(t^6+1)^{1/2} & 0 \\
 0 & -t^{1/2}(t^6+1)^{1/2} 
\end{bmatrix},
\]
for which
\[
\mu_2[A(t)+I_2K(t)]=\mu_2\left[A_{\mathrm{skew}}(t)+\mathrm{diag}\{\lambda_1,\, \lambda_2\}+\mathrm{diag}\{\gamma_1(t),\, \gamma_2(t)\}\right]
\]
\[
=\mu_2\left[\mathrm{diag}\{\lambda_1,\, \lambda_2\}+\mathrm{diag}\{\gamma_1(t),\, \gamma_2(t)\}\right]
\]
\[
=\left\{ 
\begin{array}{l}
    -1-t(t^6+1)^{1/2} \quad \mathrm{if}\ 0\leq t\leq 1\\
    -1-t^{1/2}(t^6+1)^{1/2} \quad \mathrm{if}\ t>1
\end{array} \right. 
.\]

Figure~\ref{simulation_x_1_2} depicts the simulation results for arbitrarily chosen initial state.
\end{ex} 

\

\section*{Conclusion}

By assuming knowledge of the asymptotic behavior (for $t\to\infty$) of system's modeling uncertainty $\Delta(t)$ and external disturbance represented by the nonlinear perturbing term $\omega(x,t)$, a robust adaptive state-feedback control law $u=K(t)x$ was constructed for a control system
\[
\dot x=\big[A(t)+\Delta(t)\big]x+Bu+\omega(x,t),\ t\geq t_0,
\]
to preserve the (global) asymptotic stability of the system without external disturbance, i.~e., $\omega(x,t)\equiv0.$ Surprisingly, even in the situations when the origin $x=0$ is not an equilibrium point of the perturbed system and for an asymptotically unbounded disturbances $\omega$ (Example~\ref{example_all_to_zero}), the proposed controller keeps the system asymptotically insensitive to the disturbances, in the sense of convergence of all solutions to $0$ as $t\to\infty.$


\end{document}